\documentclass[11pt, oneside]{article}   	
\usepackage{eucal}
\usepackage{graphicx}
\usepackage{color}
\usepackage{amsopn,amssymb,enumitem}
\usepackage{marginnote}
\usepackage{amsmath}
\usepackage{setspace}
\newtheorem{thm}{Theorem}[section]
\newtheorem{prop}[thm]{Proposition}

\newtheorem{lemma}[thm]{Lemma}
\newtheorem{coro}[thm]{Corollary}

\newcommand{\Sh}{\mathcal{S}}
\newcommand{\ddb}{\overline{\partial}}
\newcommand{\dd}{\partial}
\newcommand{\pa}{\partial}
\newcommand{\C}{\mathbb{C}}

\newcommand{\B}{\mathbb{B}}
\renewcommand{\H}{\mathbb{H}}

\newcommand{\R}{\mathbb{R}}

\newcommand{\Rn}{\R^n}

\newcommand{\Rnp}{\R^{n+1}}

\newcommand{\Bnp}{\B^{n+1}}

\newcommand{\Sn}{\mathbb{S}^n}

\newcommand{\St}{\mathbb{S}^2}
\newcommand{\Hnp}{\H^{n+1}}
\newcommand{\Ht}{\H^2}
\newcommand{\Hsp}{\H^3}

\newcommand{\Id}{\operatorname{Id}}

\renewcommand{\phi}{\varphi}
\renewcommand{\emptyset}{\varnothing}

\newcommand{\tr}{\operatorname{tr}}

\DeclareMathOperator{\ch}{ch}
\DeclareMathOperator{\sh}{sh}
\def\eproof{$\Box$ \medskip}

\newcommand{\Gm}{\mathcal G}
\usepackage{comment}
\numberwithin{equation}{section}
\usepackage{fancyhdr}

\title{Envelopes of Horospheres and Weingarten Surfaces in Hyperbolic 3-Space}
\author{Charles L. Epstein\footnote{This research was supported in part by an NSF Postdoctoral Fellowship.} \\ Department of Mathematics \\ Princeton University \\ Princeton, N.J.}
\date{}							

\begin{document}
\setlength{\abovedisplayskip}{0pt}
\setlength{\belowdisplayskip}{0pt}

\maketitle
\begin{abstract}We derive basic differential geometric formul{\ae} for surfaces in hyperbolic space represented as envelopes of horospheres. The dual notion of parallel hypersurfaces is also studied. The representation is applied
to prove existence and regularity theorems for Weingarten surfaces in $\mathbb{H}^3$ which satisfy
$$(1-\alpha)K = \alpha (2-H), $$ for an $\alpha<0,$ and have a specified
boundary curve at infinity. These surfaces are shown to be closely connected to
conformal mappings of domains in $\mathbb{S}^2$ into the unit disk and provide
Riemannian interpretations for some conformal invariants associated to such
mappings.

This paper was originally written in 1984, before I learned to use TeX, and was
typed by one of the secretaries in the Princeton Math Department. It was more or
less, my first original work after my dissertation. For some reason, I was not
able to get this paper published in a timely manner, and it was consigned to
what eventually became a long list of unpublished manuscripts. Some parts of
this paper appeared in an Appendix to \cite{PaPe}.

The results and
perspective in this paper have proved to be useful to a variety of
people, some of whom asked me to render the article into TeX and post
it to the arXiv. I had been seriously thinking about doing this, when
Martin Bridgeman sent a transcription of my original article into
TeX. I am extremely grateful to him for the effort he has put into
this project.

The paper is now formatted in a more or less modern AMS-article style, but for lots of
additional punctuation, a few corrections and some minor stylistic changes, the
content has been largely reproduced as it originally was.  Remarks about
the ``state-of-the-art'' in hyperbolic geometry are obviously way out of date,
as there has been enormous progress in many aspects of this still rich subject.

I am enormously grateful to Martin and the community of mathematicians
who have let me know, over the years, that this work was of some use
to them.

\end{abstract}

\section*{}
\setcounter{section}{1}
\setcounter{page}{1}

The theory of immersed surfaces in hyperbolic space is a rich and
largely unexplored subject. Recently several authors have studied surfaces
of constant mean curvature, see: \cite{Br}, \cite{DoLa}, \cite{Mi}, \cite{Uhl}. In this paper
we present a representation for hypersurfaces in $\mathbb{H}^{n+1}$ 
as graphs over the ideal boundary of $\mathbb{H}^{n+1}.$

In three dimensions three classes of surfaces are distinguished in this representation by the simplicity of their defining equations. If $k_1,k_2$ denote the principal curvatures of the immersed surface $\Sigma,$
the Gauss and Mean curvatures are:
\begin{eqnarray*}
K&=k_1 k_2-1\\
H&= k_1 + k_2.
\end{eqnarray*}
Three distinguished classes are given by the curvature conditions:
\begin{description}[noitemsep]
\item A. $K= 0,$
\item B. $H= 2,$
\item C. $(1-\alpha)K=\alpha(2-H)\ ;\alpha \in \R\setminus \{0,1\}.$
\end{description}
Surfaces of type A. have not been considered in print. They are, in Thurston's
language, flat orbifolds. Surfaces of type B. were considered by Bryant in
\cite{Br} from a point of view which is related to ours though arrived at
independently. Other results on surfaces of this type can be found in
\cite{DoG}, \cite{DoLa}. We will consider the third type of surface for
$\alpha<0.$

This paper is divided into three parts. In the first part, $\S\S$2-3, we present
the representation theory for hypersurfaces in $\H^{n+1}$ as envelopes of
horospheres, and derive the basic differential geometric formul{\ae} in this
representation. In the second part, $\S\S$4-6, refinements and extensions of the
theory in $\S\S$2-3 available for surfaces in $\H^3$ are explored. Finally, in
$\S\S$7-8 the theory developed in the first 2 parts is applied to study a
Dirichlet problem for
surfaces which satisfy:
\begin{equation}
(1-\alpha)K = \alpha(2-H)\quad\text{ with } \alpha <0. 
\end{equation}
We call such a surface an $\alpha$-Weingarten surface.

The results in the latter sections are obtained through a connection
between $\alpha$-Weingarten surfaces and conformal maps of $\St$
into $\C$. This is somewhat in the same spirit as the work in \cite{Br}. We will obtain
Riemannian interpretations for various conformal properties of such maps. 

\begin{center}
{\bf Definitions}
\end{center}

Hyperbolic $(n+1)$-space, $\Hnp$ will be represented as the interior of the unit ball, $\Bnp$ in $\Rnp$ with the metric:
$$ ds^2 = \frac{4(dx_1^2 + \ldots + dx_{n+1}^2)}{(1-r^2)^2}.$$

The ideal boundary of hyperbolic space, $\partial \Hnp$ is naturally identified with the unit sphere, $\Sn$ in $\Rnp$. The basic facts of hyperbolic geometry will be taken for granted; as a reference one can consult \cite{Th}, \cite{Be}, \cite{Sp}.

In what follows the horospheres play a central role. They are the simply
connected, complete, flat hypersurfaces in $\Hnp$.
In the ball model they are represented by Euclidean spheres internally tangent to the unit sphere.

\begin{figure}[htbp] 
   \centering
   \includegraphics[width=3in]{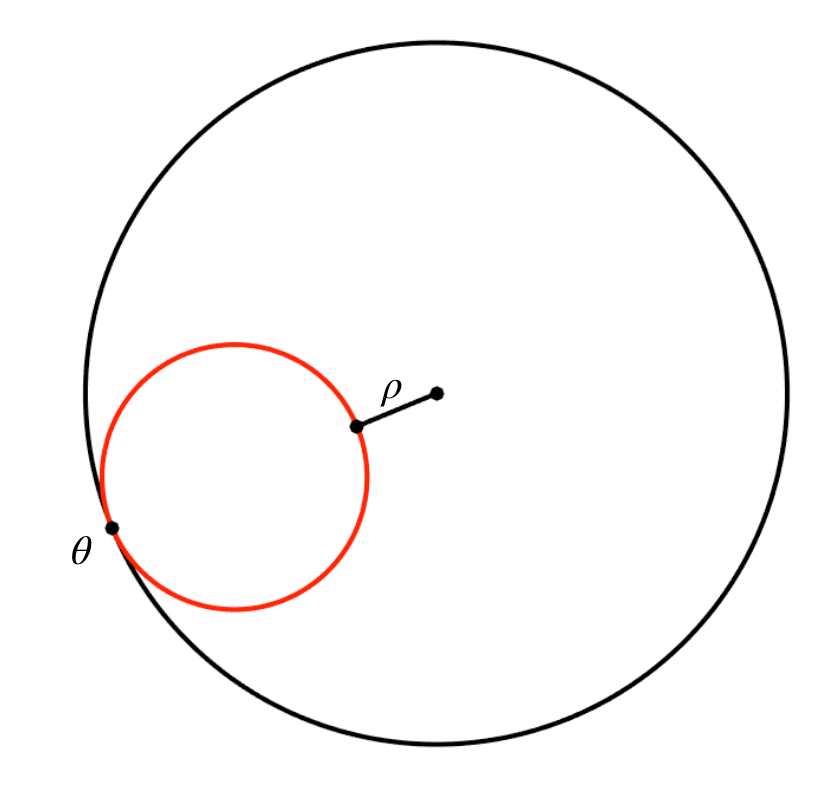} 
   \caption{The horosphere $H(\theta,\rho)$  in $\Ht$.}
   \label{fig1}
\end{figure}
\noindent
They are parameterized by $\theta$, the point of tangency with $\Sn$ 
and $\rho$, the smallest hyperbolic distance between the horosphere and the point
$(0,\dots,0)$ in  $\Bnp$.
Here $\rho$ is positive if $(0,..,0)$ is in the exterior of the horosphere and negative otherwise. We denote this
horosphere by $H(\theta, \rho)$.

Recall also that the geodesics in the ball model are the circular
arcs in $\Hnp$ that meet $\Sn$ normally. Every pair of points on $\Sn$ uniquely determines a geodesic and vice versa. We use the notation $\psi^t(p,X)$ to denote the geodesic with initial point $p\in\Hnp$ and velocity $X\in T_p\Hnp$ ; this is often denoted $\exp_p(tX)$.

If $\Sigma$ is an oriented hypersurface, smoothly embedded in $\Hnp$, then there
is a globally defined unit normal field, $N$. Using the unit ball
model, one can define a Gauss map for $\Sigma$ by:
$$\Gm_\Sigma(p) = \lim_{t\rightarrow \infty} \psi^t(p,N).$$
$\Gm_\Sigma$ has certain properties which are analogous to the Euclidean Gauss map.
These were studied independently in [Br]. Note that one could also use
$-N$ as the unit field on $\Sigma$. The mapping one obtains can be quite
different. If $A$ is an isometry of hyperbolic space, then $A\cdot \Sigma$ is also an oriented immersed surface and
$$\Gm_{A\cdot\Sigma}(Ap) = A\cdot\Gm_\Sigma(p).$$

\begin{center}
{\bf Notation}
\end{center}

Some of the notation will not be used until the end of the paper; we include it for the convenience of the reader.


\begin{eqnarray*}
 |\cdot | &-& \mbox{ Length in the Euclidean metric.}\\ 
 \langle\cdot, \cdot\rangle &-& \mbox{ The inner product on } T\Hnp\\
D &-& \mbox{ The gradient w.r.t. the round metric on }  \Sn.\\
\nabla_XY &-&\mbox{ Covariant differentiation in } \Hnp  \mbox{ w.r.t. } \langle, \rangle.\\
R(X,Y)Z &=& \nabla_X\nabla_Y Z -\nabla_Y\nabla_Z Z -\nabla_{[X,Y]}Z\\
&=& -[ \langle X,Z\rangle Y - \langle Y,Z\rangle X], \mbox{ the Riemann tensor of } \nabla_XY.\\
\Delta_S &-& \mbox{ The Laplace-Beltrami operator on } S\\
\Sigma &-& \mbox{ an immersed hypersurface}.\\
X_1,\ldots,X_n &-& \mbox{ The vector fields spanning  $T\Sigma$ defined by an immersion}.\\
N &-&\mbox{ The unit normal field of $\Sigma$}.\\
g_{ij} = \langle X_i,X_j\rangle &-& \mbox{ The induced metric on  $T\Sigma$}.\\
g^{ij} &-&\mbox{ The dual metric:} \qquad g^{ik}g_{kj} = \delta^i_j\\
\Pi_{ij} &=& \langle \nabla_{X_j}X_i, N\rangle\\ 
\Pi^i_j &=& g^{ik}\Pi_{kj}\mbox{- The second fundamental form of $\Sigma$}.\\
\end{eqnarray*}

\begin{eqnarray*}
  K(g)&-&\mbox{ The Gauss curvature of the metric $g$}.\\
D_r&-& \mbox{The disk of radius $r$ in $\C$}.\\
  \Omega &-& \mbox{ A domain in $\Sn$}\\
d\sigma^2 &-& \mbox{ The curvature +1 metric on $\Sn$}\\
z&-&\mbox{ A local conformal parameter on } \Sn\ . \\
'&-&\mbox{Differentiation w.r.t. } z\ .\\
\gamma&-&d\sigma^2 = \gamma^2|dz|^2\ .\\
\rho&-&\mbox{ A function defined on a domain in } \Sn\ . \\ 
\rho_\Omega&-&\mbox{ The hyperbolic metric on } \Omega\subset\St \mbox{ is } e^{2\rho_\Omega}d\sigma^2 \  . \\
\Sigma(\rho)&-&\mbox{ The surface generated by } \rho\ .\\
\Sigma_t&-&\mbox{ The parallel surface at distance $t$ from
  $\Sigma$}\ . \\
f_\Omega&-&\mbox{ The conformal map from $\Omega$ to $D_1$} \ . \\
\left.\begin{matrix} \dd = 1/2(\dd_x-i\dd_y)\\
\ddb  = 1/2(\dd_x+i\dd_y)\end{matrix}\right\}&-&\mbox{The complex derivatives w.r.t. a conformal parameter $z=xd+iy$}\ .\\
\Sh_f(z)&=&\left[\left(\frac{f''}{f'}\right)'-\frac{1}{2}\left(\frac{f''}{f'}\right)^2\right](z)\\
&-&\mbox{ The Schwarzian derivative of $f$}\\ 
\mu_\Omega&-&\mu_\Omega|dz|^2 \mbox{ the hyperbolic metric on
  $\Omega$ w.r.t. the conformal parameter $z$}\ . \\
& &\mu_\Omega = 4|f'_\Omega|^2(1-|f_\Omega|^2)^{-2}\\
\mbox{Note}&-&\mu_\Omega |dz|^2 = e^{2\rho_\Omega}d\sigma^2 \mbox{ if $d\sigma^2$ is represented w.r.t. $z$}\ . 
\end{eqnarray*}

\begin{center}
{\bf Acknowledgments}
\end{center}

I would like to thank Bill Thurston who suggested the representation
of hypersurfaces explored and exploited in this paper and for his
continued interest and advice. I would also like to thank Robert
Bryant for allowing me to read his beautiful manuscript \cite{Br} on
surfaces of mean curvature 2 and for the conversations we had on this
and related subjects.

\section{Hyperbolic Graphs and Parallel Surfaces:}
\indent Let $\Omega$ be a domain on $\Sn$ and $\rho(\theta)$ a
differentiable function defined in $\Omega$. Recall that if $\mathfrak
A$ is a family of hypersurfaces in $\Rnp$ then the envelope of
$\mathfrak A$ is a hypersurface $\Sigma$ which is everywhere tangent
to hypersurfaces in $\mathfrak A$. This notion is independent of the
metric on $\Rnp$ as it is simply a statement about identity of tangent
spaces. A function $\rho(\theta)$  defines a surface $\Hnp$ by the prescription:
\begin{multline*}
\Sigma(\rho) =\\ \mbox{ Outer envelope of the family of }
\mbox{horospheres: } \{H(\theta,\rho(\theta)) : \theta \in \Omega\}.
\end{multline*}
Outer envelope refers to the component of the envelope that is not the unit sphere, which is
obviously an envelope for any smooth family of horospheres.

The outer envelope of a smooth family of horospheres may fail to exist in that
$\Sigma(\rho)$ may fail to be smooth. Notwithstanding, we will derive a formula
for $\Sigma(\rho)$ which gives the usual envelope wherever it exists. At other
points it defines a continuous mapping from $\Omega$ into $\Hnp$.

We need an explicit formula for the horosphere $H(\theta,\rho)$. Let
\begin{equation}
r = \frac{e^\rho-1}{e^{\rho}+1}
\label{rdef}
\end{equation}
$$H(\theta,\rho) = \left\{ \frac{1+r}{2}X(\theta)+ \frac{1-r}{2}Y\ : \ Y \in \Sn\rightarrow \Rnp\right\}.$$
$X(\theta)$ is the point $\theta$ on the unit sphere in  $\Rnp$. We will also use $\theta$ as the
notation for a coordinate.

To derive the formula for the envelope we let $Y\in \Sn$ be represented
parametrically as $Y(\alpha)$ where $\alpha=(\alpha^1,\ldots,\alpha^n)$ ranges over an open set
in $\Rn$; $\theta = (\theta_1,\ldots,\theta_n)$ also ranges over an open set in $\Rn$. 
Let

$$R(\theta,\alpha) = \frac{1+r(\theta)}{2} X(\theta) + \frac{1-r(\theta)}{2} Y(\alpha). $$
Formally to solve for the envelope we should solve for $\alpha(\theta)$. In fact
this will be unnecessary as we can solve for $Y(\alpha(\theta))$ directly. For the
remainder of the derivation we will assume that the envelope exists, i.e.,
that one could solve for $\alpha(\theta)$
and it would be differentiable. 

The tangent space to $\Sigma(\rho)$ at $R(\theta,\alpha(\theta))$ is spanned by:
$$2R_{\theta_i} = r_{\theta_i} X + (1+r)X_{\theta_i} -r_{\theta_i}Y + (1-r)\sum_{j=1}^n Y_{\alpha^j} \alpha^j_{\theta_i},\qquad i=1,\ldots,n\ . $$
The tangent space to $H(\theta,\rho(\theta))$ is spanned by:
$$2R_{\alpha^j} = (1-r)Y_{\alpha^j},\qquad j=1,\ldots n\ .$$
The conditions defining the envelope are
\begin{equation}
\mbox{Span} \{ R_{\theta_i}(\theta,\alpha(\theta)) \} = \mbox{ Span}\{ R_{\alpha^i}(\theta,\alpha(\theta))\}\ .
\label{Rspan}
\end{equation}
This is true at a regular point if and only if:
$$R_{\theta_i}(\theta,\alpha(\theta)) \perp R_{\alpha^1}\times\ldots\times R_{\alpha^n}$$
where $\perp$ and $\times$ are with respect to the Euclidean inner product, which we denote by $X \cdot Y$.

An elementary calculation shows that
$$R_{\alpha^1}\times\ldots\times R_{\alpha^n} = \lambda Y$$
for some $\lambda \neq 0$. We rewrite the conditions, \eqref{Rspan} as:
\begin{equation}
r_{\theta_i}X\cdot Y + (1+r)X_{\theta_i}\cdot Y = r_{\theta_i}\qquad 1=1,\ldots,n. 
\label{2.2} \end{equation}
These follow as $Y\cdot Y = 1$ and $Y\cdot Y_{\alpha^j} = 0$. This is an inhomogeneous system
of linear equations for $Y(\theta)$. If we compute the Gramian matrix associated to the system we obtain:
$$G = (1+r)^2\Id + r_\theta\otimes r^t_\theta.$$
Since $G$ is always of rank $n$, the kernel of \eqref{2.2}  is one dimensional.
It is generated by the vector:
$$Z = X- (r+1)^{-1}\sum_{i=1}^n r_{\theta_i}X_{\theta_i}(X_{\theta_i}\cdot X_{\theta_i})^{-1}.$$
A particular solution to \eqref{2.2} is $Y_p = X$, so the general solution is 
$$X+\mu Z\qquad \mu \in \R.$$

The condition that we use to determine $\mu$ is $Y\cdot Y= 1$.  
The possible solutions are:
\begin{eqnarray}
Y &=& X\qquad \mbox{and} \nonumber \\
Y &=& \frac{|Dr|^2-(1+r)^2}{|Dr|^2+(1+r)^2}X + \frac{2(1+r)Dr}{|Dr|^2+(1+r)^2}
\end{eqnarray}
where
\begin{eqnarray*}
Dr &=& \sum r_{\theta_i}X_{\theta_i}(X_{\theta_i}\cdot X_{\theta_i})^{-1} \mbox{ and}\\
|Dr|^2 &=& Dr\cdot Dr.
\end{eqnarray*}
$D$ is the gradient on $\Sn$ with respect to the round metric. $Y = X$ is clearly the inner envelope and thus using \ref{rdef} we obtain a formula for $R_\rho$:
\begin{equation}
R_\rho(\theta) = \frac{|D\rho|^2+(e^{2\rho}-1)}{|D\rho|^2+(e^{\rho}+1)^2}X(\theta) + \frac{2D\rho}{|D\rho|^2+(e^{\rho}+1)^2}.
\label{2.4}
\end{equation}
We will use $R_\rho(\theta)$ to denote the parametric representation of the
mapping from $\Omega$ into $\Hnp$ 
defined by $\rho$. The parametrization is determined by the parametrization $X(\theta)$ of the unit sphere.

From the derivation it is clear that $R_\rho(\theta)$ coincides with the envelope of $H(\theta, \rho(\theta))$ whenever it is possible to solve for $\alpha(\theta)$. For
under this assumption there is a unique point on each horosphere which lies in the envelope. We will use $\Sigma(\rho)$  to denote the hypersurface generated
by $\rho$. As we will see in $\S$5, if $\rho$ is twice differentiable then $\alpha(\theta)$
exists wherever $\Sigma(\rho)$ is smooth. Clearly $R_\rho(\theta)$ is continuous if $\rho(\theta)$ is continuously differentiable and thus the envelope is connected if $\Omega$ is connected.

If $\rho \in C^k(\Omega)$, then $\Sigma(\rho)$ will be $C^{k-1}$-surface wherever $R_\rho$ is
an immersion. An interesting feature of this representation is that the principal curvatures of $\Sigma(\rho)$ are expressions involving $\rho$ and its first and second derivatives. As the formula for $R_\rho$ involves $\rho$ and its first derivative, we would expect the principal curvatures to depend on the
third derivatives of $\rho$  as well.

Closely connected to the construction of hypersurfaces as envelopes of horospheres is the family of hypersurfaces parallel to a given hypersurface. Let $N$ denote the outward unit normal field on $\Sigma$. Then the parallel surface at distance $t$ is defined to be:
$$\Sigma_t = \{ \psi^t(p,N)\ :\ p \in \Sigma\}\ .$$
The dynamical system $\Sigma \rightarrow \Sigma_t$ will be referred to as the parallel flow.
The connection between parallel surfaces and envelopes of horospheres is contained in the next theorem. This fact was observed by Thurston:

\begin{thm}
The parallel surface at distance $t$ from $\Sigma$ is the surface
generated by the function $\rho_t=\rho + t$ ; more succinctly:
$$\Sigma_t(\rho) = \Sigma(\rho_t).$$
\end{thm}

\noindent\underline{Remark:} $\Sigma(\rho)$ need not be smooth; at a non-smooth point there is no normal
vector. However, the representation $R_\rho$ allows us to define a subspace of $T\Hnp$ at non-smooth points composed of ``normal" vectors to $\Sigma(\rho)$.
If $\mathcal N_P$ is this family of vectors then
$$S_t = \{\psi^t(p,N)\ : \  N \in \mathcal N_p\}$$
will be a smooth hypersurface for small non-zero $t$. Fix such a $t = t_0$. Then $p \in \Sigma(\rho)$  will be a focal point of the surface $S_{t_0}$ (see $\S$3). At smooth points on $\Sigma(\rho)$ the statement of the theorem is rigorously correct. We will make further remarks after the proof.

\noindent\underline{Proof of Theorem 2.1:}
Consider the formula for $R_{\rho+t}(\theta)$ :
$$R_{\rho+t}(\theta) = \frac{|D\rho|^2+(e^{2(\rho+t)}-1)}{|D\rho|^2+(e^{\rho+t}+1)^2}X(\theta) + \frac{2D\rho}{|D\rho|^2+(e^{\rho+t}+1)^2}.
$$
Elementary but tedious calculations show that for fixed $\theta$:
\begin{enumerate}
\item $\left|R_{\rho+t}(\theta) - (X(\theta)+\frac{D\rho(\theta)}{|D\rho(\theta)|^2}|)\right| \mbox{ is constant.}$

$|\cdot |$ is hyperbolic distance.
\item $$\lim_{t\rightarrow \infty} e^t \frac{dR_{\rho+t}(\theta)}{dt} = \lim_{t\rightarrow \infty} 2R_{\rho+t}(\theta) = 2X(\theta)$$
and
\begin{eqnarray*}
\lim_{t\rightarrow -\infty} e^{-t} \frac{dR_{\rho+t}(\theta)}{dt} &=& \lim_{t\rightarrow -\infty} -2R_{\rho+t}(\theta)\\
&=& -2 \left[ \frac{|D\rho|^2-1}{|D\rho|^2+1}X(\theta) + \frac{2D\rho}{|D\rho|^2+1}\right].
\end{eqnarray*}
The limits are taken in the Euclidean topology of $\Rnp$.
\item $$\left\langle \frac{dR_{\rho+t}(\theta)}{dt},\frac{dR_{\rho+t}(\theta)}{dt}\right\rangle = 1 \mbox{ for all $t$}.$$
\end{enumerate}
From 1), 2) and 3) it follows that $R_{\rho+t}(\theta)$ is a unit speed parametrization of a hyperbolic geodesic normal to the family of horospheres $\{H(\theta,\rho(\theta)+t)\ :\ t\in \R\}$. If $\Sigma(\rho)$ is smooth at $p$ then $\Sigma(\rho)$ is tangent to $H(\theta,\rho(\theta))$ at $p$ and  thus the unit normal vector to $\Sigma(\rho)$ at $p$ is:
$$N = \left.\frac{dR_{\rho+t}(\theta)}{dt}(\theta)\right|_{t=0}.$$

From the definition of a parallel surface, it is clear that if $p= R_{\rho}(\theta)$ is a smooth point on $\Sigma(\rho)$ then:
$$R_{\rho+t}(\theta) = \psi^t(p,N).$$
Thus the theorem follows wherever $\Sigma(\rho)$ is smooth.

If $\Sigma(\rho)$ is not smooth at $R_\rho(\theta)$ then $R_{\rho+t}(\theta)$ is
a point a distance $t$ from $\Sigma(\rho)$. If $\Sigma(\rho+t)$ is smooth at
this point, then the normal vector to $\Sigma(\rho+t)$ at $R_{\rho+t}(\theta)$
is
$$\widetilde{N} = \left.\frac{dR_{\rho+t+s}(\theta)}{dt}(\theta)\right|_{s=0}.$$
$R_\rho(\theta)$ is the point on $\Sigma_{-t}(\rho+t)$ given by:
$$\psi^{-t}(p,\widetilde{N}).$$
$R_\rho(\theta)$ is in the focal set of $\Sigma(\rho+t)$. We will discuss this further  $\S$3.

\noindent\underline{Remark}: From the argument in the proof it is clear that the family of
hypersurfaces: $\{\Sigma_t(\rho)\ :\ t\in\R\}$  agrees insofar as it is defined with the family of hypersurfaces:
$$\{\Sigma(\rho+t)\ :\ t\in\R\}\ .$$
Moreover, the second family provides a consistent completion for the first
family.

One could eliminate this discussion altogether by letting $\rho$ define a mapping into the unit tangent bundle, $T_1\Hnp$ by
$$\theta \rightarrow \left(R_{\rho}(\theta), \left.\frac{dR_{\rho+t}(\theta)}{dt}(\theta)\right|_{t=0}\right) \in T\Hnp|_{\Sigma(\rho)}\ .$$
The parallel flow is defined in an obvious way on $T_1\Hnp$ and is
everywhere smooth; it is the projection into $\Hnp$ that produces the
singularities. In this context the statement:
$$T\Sigma_t(\rho) = T\Sigma(\rho+t)$$
is rigorously correct at every point. 
\section{Equations of Motion for Parallel Surfaces}
In this section we will derive equations that describe the evolution
of the fundamental geometric quantities of an immersed hypersurface under
the parallel flow. In particular, we will derive equations of motion for
the principal curvatures. The rather startling fact is that these equations
are decoupled; consequently, the principal curvatures evolve independently
of one another. This allows a good qualitative understanding of the
formation of singularities under the parallel flow,

We will not employ the representation discussed in $\S$2, rather we
will consider hypersurfaces as given locally by embeddings of open sets
$U \subset \Rn$: 
$$i: U \rightarrow \Hnp.$$
If $(x_1,\ldots,x_n)$ are coordinates on $U$ then the vector fields:
$$X_j = i_*\partial_{x_j}\qquad j= 1,\ldots,n$$
span the tangent space of $i(U)$. $N$ will denote the outward unit normal
vector field, The first fundamental form is:
$$g_{ij} = \langle X_i,X_j\rangle$$
the second fundamental form is:
$$\Pi^i_j = g^{ik}\Pi_{kj},$$
where
$$\Pi_{kj} = \langle \nabla_{X_j} X_i ,N\rangle.$$
$\{X_i,\, \ i=1,\ldots,n\}$ and $N$
will also be used to denote respectively the coordinate and normal vector fields
of the parallel hypersurfaces. The extensions of the $X_i$ to the parallel
hypersurfaces are defined by the immersion:

$$\psi^{t}(i(x),N_{i(x)})\qquad\qquad x\in U.$$
As the tangent vector field to a geodesic,
$$N = \frac{d\psi^t}{dt}(i(x),N_{i(x)})$$
satisfies the geodesic equation
$$\nabla_N N = 0\ .$$

All the commutators, $[X_i,X_j],\  i,j=i,\ldots,n$ and $[X_i,N], \ i=1,\ldots,n$
vanish. From this we easily conclude that
$$\langle X_i,N\rangle = 0$$
for all $t$.

The tangential vector fields, $\{X_1,\ldots,X_n\}$ satisfy the Jacobi equation:
\begin{equation}
\nabla^2_NX_i + R(N,X_i)N = 0\ .
\end{equation}
The quantities on the parallel hypersurfaces corresponding to those
introduced above will also be denoted by $g_{ij}, \Pi_{ij}$ and
$\Pi^i_j$ respectively, Whenever we want to emphasize the dependence
on a particular variable, we will write a quantity as explicitly
depending on that variable. For instance $g_{ij}(t)$; $\frac{d}{dt}$
denotes the action of the vector field N.

\begin{thm} The first and second fundamental forms satisfy the following equations under the parallel flow:
\begin{enumerate}[label=\alph*)]
\item $$\frac{dg_{ij}}{dt} = -2\Pi_{ij} $$
\item 
\begin{equation}
\frac{d\Pi^i_j}{dt} = \Pi^i_l\Pi^l_j-\delta^i_j
\end{equation}
\item 
$$\frac{d\Pi_{ij}}{dt} = 4\Pi_{ij}.$$
\end{enumerate}
\end{thm}

\noindent\underline{Remark:} The second fundamental form is $\Pi^i_j$ and it evolves
independently of $g_{ij}$ . Note however, that (a) and (c) are linear
while (b) is non-linear, For this reason it is sometimes advantageous
to use (a) and (c) instead of (b).

\noindent\underline{Proof of Theorem 3.1}: The proofs are elementary calculations using the
properties of the Levi-Civita connection;
\begin{eqnarray*}
a)\  \frac{dg_{ij}}{dt} &=& N\langle X_i,X_j\rangle\\
&=& \langle \nabla_NX_i,X_j\rangle+\langle X_i,\nabla_NX_j\rangle\\
&=& \langle \nabla_{X_i}N,X_j\rangle+\langle X_i,\nabla_{X_j}N\rangle,\\
\mbox{as } [X_i,N] = 0 &\text{ for }& i=1,\ldots,n. \mbox{ This equals}\\
&=& -\left[ \langle N,\nabla_{X_i}X_j\rangle+\langle N,\nabla_{X_j}X_i\rangle\right]\\
&=& -2\Pi_{ij}.
\end{eqnarray*}
b) To prove part (b) it is more convenient to begin with an orthonormal frame
adapted to the hypersurfaces, Let $\Sigma$ denote the initial surface and
$\{Y_1,\ldots, Y_n\}$ denote an orthonormal frame field tangent to $\Sigma$. Extend this
frame to $\Sigma_t$ by parallel translation along the normal geodesics to $\Sigma$,

Denote the extended vector fields by $\{Y_1,\ldots, Y_n\}$ as well. They are
orthonormal and satisfy:
\begin{equation}
\begin{matrix}
\nabla_NY_i &=& 0\\
\langle N, Y_i\rangle &=& 0 
\end{matrix}\qquad \text{ for } i=1,\ldots,n.
\end{equation}
In terms of such a frame:
\begin{equation}
\overline\Pi^i_j = \langle \nabla_{Y_i}Y_j, N\rangle.
\end{equation}
We differentiate (3.4) to obtain:
$$
\frac{d\overline\Pi^i_j }{dt}= \langle \nabla_N\nabla_{Y_i}Y_j, N\rangle.
$$
Using the definition of the Riemann tensor, this can be rewritten as
\begin{equation}\begin{split}
\frac{d\overline\Pi^i_j }{dt} & = \langle \nabla_{Y_i}\nabla_{N}Y_j, N\rangle+\langle \nabla_{[N,Y_i]}Y_j, N\rangle\\
&+ \langle R(N,Y_i)Y_j, N\rangle\\
&= \langle \nabla_{[N,Y_i]}Y_j, N\rangle-\delta^i_j.
\end{split}\end{equation}
We have used (3.3) and the formula for the Riemann tensor of $\Hnp$
in the introduction. The connection is symmetric thus;
$$[X,Y] = \nabla_XY - \nabla_YX.$$
Hence it follows from (3.3) that
$$[N,Y_i] = - \nabla_{Y_i}N.$$
$\nabla_YN$ is tangent to $\Sigma_t$ and so can be expressed in terms of the frame $\{Y_1,\ldots, Y_n\}$, by 
\begin{equation}
[N,Y_i] = \overline{\Pi}^i_kY_k.
\end{equation}
Using (3.6) in (3.5) we obtain:
$$\frac{d\overline\Pi^i_j}{dt} = \overline\Pi^i_l\overline\Pi^l_j-\delta^i_j$$
as asserted. To complete the derivation of b we introduce matrix notation:
\begin{eqnarray*}
\overline P &= \overline{\Pi}^i_j\\
 P &= \Pi^i_j\\
Q &=\Pi_{ij}\\
G &= g_{ij}.
\end{eqnarray*}
In this notation 
\begin{equation}
P=G^{-1}Q.
\end{equation}
As $\{Y_1,\ldots, Y_n\}$ is a basis for each tangent space $T_p\Sigma_t$ ,  there is a
matrix $A_{ij}$, such that
$$ X_i = \sum_j A_{ij}Y_j.\qquad\qquad (*) $$
We can express $G, Q$ and $P$ in terms of $A$ and $\overline P$:
\begin{eqnarray*}
G&=& AA^t\\
 Q &=& A\overline P A^t\\
P &=& (A^t)^{-1}\overline P A^t \ .
\end{eqnarray*}
We differentiate $(*)$ with respect to time to obtain:
$$\nabla_NX_i = \sum \frac{dA_{ij}}{dt} Y_j\ .$$
By taking the inner product with $Y_k$ and applying $(*)$ we obtain:
$$\frac{dA}{dt} =-A\overline P\ .$$
Using these formul{\ae} and the equation satisfied by $\overline P$ one easily derives 
$$\frac{dP}{dt} = P^2 - I .$$
which completes the proof of (b).
To derive c we
use the matrix notation introduced above:
\begin{eqnarray*}
\frac{dP}{dt} &= -G^{-1}\frac{dG}{dt}G^{-1}Q+G^{-1} \frac{dQ}{dt}\\ &= 2P^2 +G^{-1}\frac{dQ}{dt}.
\end{eqnarray*}
Recall that
$$\frac{dP}{dt} = P^2-I$$
thus
\begin{equation}
\frac{dQ}{dt} -(PQ + G).
\end{equation}
We differentiate (3.8) to obtain:
$$\frac{d^2Q}{dt^2}=\frac{dQ}{dt}P+Q\frac{dP}{dt}+\frac{dG}{dt}.$$
Using (3.8) and (3.2)(a) and (b) we obtain:
$$\frac{d^2Q}{dt^2}=4Q.$$\eproof

Using (3.2) (a) and (c) it is possible to write an explicit formul{\ae} for
$g_{ij}$, and $\Pi_{ij}$ :
\begin{eqnarray*}
\Pi_{ij} &=& e^{2t}\Pi_{ij}^++e^{-2t}\Pi_{ij}^-\\
g_{ij} &=& \Gamma_{ij}-e^{2t}\Pi_{ij}^+ +e^{-2t}\Pi_{ij}^-.
\end{eqnarray*}
$\Pi_{ij}^\pm$ and $ \Gamma_{ij}$ are functions of position on the initial
hypersurface, The principal curvatures $\{k_1,\ldots,k_n\}$ are the eigenvalues
of $\Pi^i_j$ , they are real as $\Pi^i_j$ is symmetric in an orthonormal frame.
The eigenvectors of $\Pi^i_j$ are the principal directions. Corollary 3.2
follows easily from equation (3.2)(b):

\begin{coro} a) The principal curvatures satisfy the equations:
\begin{equation}
\frac{dk_i}{dt} = k_i^2-1 \qquad i=1,\ldots,n.
\end{equation}
b) The principal directions on $\Sigma_t$  are the images under the parallel flow of the principal directions on $\Sigma$.
\end{coro}
\noindent\underline{Proof of Corollary 3.2}: As $\Pi^i_j$ is a $1 \choose 1$-tensor, it is diagonal if and only if it is represented by a frame composed of eigenvectors, that is principal directions, In this case:
$$\Pi^i_j = k_i \delta^i_j\ .$$
The inhomogeneity in (3.2)(b) is a diagonal matrix. Thus a solution
with diagonal initial data is diagonal for all time. The diagonal entries satisfy the equation
$$
\frac{dk_i}{dt} = k_i^2-1 \qquad i=1,\ldots, n,
$$
proving the assertion in part (a).

To prove part (b) we observe that the coordinates at a point $p$ 
on $\Sigma$ can be chosen so that the initial coordinate vector fields $\{X_1,\ldots, X_n\}$ are the principal directions at $p$. Thus $\Pi^i_j(0,p)$  will be diagonal and so it will be diagonal for all time. $\Pi^i_j(t,p)$ is represented with respect to the extended coordinate vector fields $\{X_1(t,p),\ldots, X_n(t,p)\}$.
From the remarks in the proof of part (a) it follows that $\{X_1(t,p),\ldots, X_n(t,p)\}$ must be a frame of principal directions for  $\Sigma_t$ at $\psi^t(p,N)$.

If the principal curvatures are distinct at a point $p\in \Sigma$ then we can construct smooth vector fields $\{P_1,\ldots, P_n\}$ composed of principal directions in a neighborhood of $p$. The integral curves to these vector
fields are called the lines of curvature of $\Sigma$. Part (b) of the Corollary could be rephrased: The lines of curvature on $\Sigma_t$ are the images under the parallel flow of the lines of curvature on $\Sigma$.
Using equation (3.9) we can derive a formula for $k_i(t)$ :
\begin{equation}
k_i(t) = \frac{k_i(0)\ch t-\sh t}{-k_i(0)\sh t+\ch t}.
\end{equation}
From this it is evident that $k_i(t)$ will be infinite at some time
if and only if:
$$|k_i(0)| > 1.$$

This is a reflection of the non-linear evolution equation governing
the second fundamental form. On the other hand (3.2)(a) and (c) are linear
equations and therefore $g_{ij}$ and $\Pi_{ij}$ are smooth and bounded
for all finite times. For some $t$ $g_{ij}(t)$ may fail to be invertible, at such
points $\Pi^i_j(t)$ will be unbounded. To discuss this further we require
a definition:

\noindent\underline{Definition}: Let $\Sigma$ be a smooth immersed surface with unit normal field $N$. The set of points:
$$\mathfrak F = \{\psi^t(p,N)\ :\ p\in\Sigma, t \in \R, \mbox{ and } \det
g_{ij}(t,p) = 0\}\subset \Hnp$$
is called the focal manifold of $\Sigma$.  The forward focal manifold,
$\mathfrak F_+$ is the subset of $\mathfrak F $ with $t> 0$.  The backward focal
manifold, $\mathfrak F_-$ is the subset of $\mathfrak F $ with $t\leq 0$.

Let $i: U \rightarrow \Hnp$  be an immersion defining $\Sigma$ locally. As is evident from the definition, the composed mapping
$$i_t(x) = \psi^t(i(x),N_{i(x)})$$
fails to be an immersion if $\psi^t(i(x),N_{i(x)}) \in \mathfrak F$. We will show that the surface $\Sigma$ actually is singular on the focal set:

\begin{prop} If $\mathfrak F\cap \Sigma_t$ an is not empty then $\Sigma_t$ is not a smooth hypersurface. The singular locus is the image under the parallel flow of the lines of curvature on $\Sigma$ where
$$ k_i = \frac{\ch t}{\sh t} \mbox{ for some $i$}.$$
\end{prop}

\noindent\underline{Proof of Proposition 3.3}: We need to show that whenever $\det g_{ij}(t)$ is zero some $k_i(t)$ is infinite.

Fix a point p on $\Sigma$ and introduce coordinates on $\Sigma$ so that the
coordinate vector fields $\{X_1,\ldots, X_n\}$ are the principal directions at $p$. As usual $X_i$ will also denote the extended coordinate vector fields. An easy calculation shows that at $p$ :
\begin{equation}\begin{split}
\left.\nabla_N X_i\right|_{t=0} &= -k_iX_i\\
\langle X_i,N\rangle &= 0.
\end{split}\end{equation}

Recall that the vector fields $\{X_1,\ldots, X_n\}$ are solutions of the Jacobi equation, (3.1). As $X_i(0,p)$ satisfies the initial conditions in (3.11), it can be expressed as:
\begin{equation}
X_i(t,p) = \frac 12[ (1-k_i(0))e^t+(k_i(0)+1)e ^{-t}] \tilde X_i(t)\end{equation}
where $\tilde X_i(0) = X_i(0)$ and
$$\nabla_N\tilde X_i = 0.$$
Thus we see that $\langle X_i, X_i\rangle = 0$ if and only if:
$$e^{2t} = \frac{k_i(0) +1}{k_i(0) -1}.$$

Solving for $k_i(0)$ we obtain:
$$ k_i(0) = \frac{\ch t}{\sh t} $$
as asserted. From formula (3.10) it is evident that $\langle X_i(t,p),X_i(t,p) \rangle = 0$ if and only if $k_i(t) = \infty$. As $\{X_1(t,p),\ldots, X_n(t,p)\}$ are a frame of principal directions, they are orthogonal; hence,
$$\det g_{ij}(t,p) = \prod_{i=1}^n \langle X_i(t,p),X_i(t,p) \rangle.$$
The proposition follows easily from this and Corollary 3.2. \eproof

If  $\mathfrak F_+$ is empty we say that $\Sigma$ is forward convex and backward convex if $\mathfrak F_-$ is empty. Evidently  $\Sigma$ is forward convex if and only if:
\begin{equation}\begin{split}
k_i &\leq 1\qquad\qquad  i= l,\ldots,n \\
\mbox{and backward convex } &  \mbox { if and only if}\\
k_i  &\geq -1 \qquad\qquad i=1,\ldots,n.
\end{split}\end{equation}

A remarkable feature of hyperbolic space is the existence of surfaces which are both forward and backward convex. This occurs if and only if
\begin{equation}
|k_i| \leq 1\qquad  i= l,\ldots,n.
\end{equation}

\noindent\underline{Note}: These notions of convexity are different from geodesic convexity.

We close this section with a theorem on hypersurfaces which satisfy (3.14). The proof was suggested by Thurston.
\begin{thm}
 If $\Sigma$ is a complete hypersurface in $\Hnp$ whose principal
curvatures satisfy (3.14), then $\Sigma$ has no self intersections.
\end{thm}

\noindent\underline{Proof of Theorem 3.4}: If this were not the case the completeness of $\Sigma$ would imply the existence of a non-trivial geodesic loop. That is an arc-length parametrized curve $c(t)$ on $\Sigma$ such that
\begin{equation}
\nabla_T T = K N
\end{equation}
where $T= \dot{c}(t)$ and $N$ is the normal vector to $\Sigma$ at $c(t)$; $K$ is the geodesic curvature of $c(t)$ as a curve in $\Hnp$. Moreover
$$c(0) = c(t_0)$$
for some $t_0 > 0$. Taking the inner product of (3.15) with $N$ 
we obtain:
$$\langle \nabla_T T,N\rangle = K$$
or 
$$\Pi(T,T) = K;$$
$\Pi(\cdot,\cdot)$ is the second fundamental form of $\Sigma$. The principal curvatures are the eigenvalues of $\Pi$; the Courant min-max principle implies that:
\begin{eqnarray*}
|\Pi(T,T)| &\leq& \max |k_i|\\
&\leq& 1\ .
\end{eqnarray*}
But this is impossible as we will show in Lemma 3.5. The proof of the theorem will be complete once we've proved:

\begin{lemma} Let $c(t)$ be an arc-length parametrized curve in $\Hnp$ with geodesic curvature $K$. Let $p$ denote a fixed point in $\Hnp$.   If $|K| \leq 1$ then
$$f(t) = \cosh(d(p,c(t)))$$
is a convex function. ($d(\cdot,\cdot)$ is hyperbolic distance).
\end{lemma}
An immediate corollary is:

\begin{coro} If $c(t)$ is as described in Lemma 3.5, then $c(t)$
has no self intersections.
\end{coro}

\noindent\underline{Proof of Lemma 3.5}: First we note that in the upper half space model $d(\cdot,\cdot)$ is the inverse hyperbolic cosine of a smooth function. Thus $f(t)$ is as smooth as $c(t)$. Let $T= \dot{c}(t)$ be the tangent vector. Let $p$ be the center of a geodesic normal coordinate system, $(r,\omega)$ for $\Hnp$. In such a coordinate system:
$$ds^2 = dr^2 + \sh^2 r d\sigma^2;$$
$d\sigma$  is the line element on the $n-$sphere of curvature $+1$. Let
\begin{eqnarray*}
c(t) &=& (r(t), \omega(t))\\
T &=& \dot{r}\partial_r + \dot{\omega}_i X_i,
\end{eqnarray*}
where $\{X_1,\ldots,X_n\}$  is an orthogonal frame of coordinate vector fields on $\Sn$.
\begin{eqnarray*}
\langle X_i, X_j\rangle &=& g_{ij}\delta^i_j\\
\langle X_i, \partial_r\rangle &=& 0 \qquad i = 1,\ldots, n\ .
\end{eqnarray*}
As c(t) is parametrized by arc-length,
\begin{equation}
\dot{r}^2 + \sh^2 r|\dot\omega|^2  = 1.
\end{equation}
The first Frenet-Serret equation reads:
\begin{equation}
\nabla_T T = KN.
\end{equation}
$N$ is a unit vector field along $c(t)$. We will need only the $\partial_r$
component of this equation. As $\Sn$ is totally umbilic, it is easy to see that
\begin{equation}
\langle\nabla_{X_i} X_j,\partial_r\rangle = 0 \mbox{ if } i\neq j.
\end{equation}
From the form of the metric it follows that:
\begin{equation}\begin{split}
\langle\nabla_{X_i} X_i,\partial_r\rangle & =  -g_{ii}\sh r \ch r\\
\langle\nabla_{\partial_r} X_i,\partial_r\rangle & =  0\\
\nabla_{\partial_r} \partial_r  &= 0.
\end{split}\end{equation}
Taking the inner product of (3.17) with $\partial_r$ and using (3.18) and (3.19), we obtain:
$$\ddot{r} = |\dot\omega|^2\sh r\ch r + K\langle\partial_r,N\rangle.$$
Thus
\begin{eqnarray*}
\frac{d^2f}{dt^2} &=& \ddot{r}\sh r + \dot{r}^2 \ch r\\
&=& \ch r (\dot{r}^2+ |\dot\omega|^2\sh^2 r) + K\langle\partial_r,N\rangle\sh r\\
&\geq& \ch r - |K| \sh r\\
&>& 0.
\end{eqnarray*}
whenever $|K| <1$. The inequality above follows from (3.16) and the fact that $\partial_r$ and $N$ are unit vectors. \eproof

\noindent\underline{Remark}: As a final remark we note that the conditions:
$$|k_i| \leq 1 \qquad i= 1,\ldots,n$$
are preserved under the flow defined on $(k_1,\ldots,k_n)$ by the system of equations (3.9). This is because the locus of points where
$$|k_i| = 1 \qquad i= 1,\ldots,n$$
consists of critical points. Thus, no trajectory crosses
any of the hyperplanes defined by
$$k_i = \pm 1 \qquad i= 1,\ldots,n\ .$$
\section{Surfaces in $\Hsp$} 
In this section and for the remainder of the paper we will consider the classical case of two dimensional surfaces in $\Hsp$. Recall that the Gauss and mean curvatures are expressed in terms of the principal curvatures by:
\begin{equation}\begin{split}
K &= k_1k_2 -1 \\
H &= k_1 + k_2\  .
\end{split}\end{equation}
Let $g$ denote $\det g_{ij}$. The area form with respect to the
parameters $(x,y)$ is
$$dA = \sqrt{g} \ dx\wedge dy\ .$$
In this section we derive equations of motion for $K, H$ and $g$.

\begin{prop} Under the parallel flow the functions $K, H$ and $g$
satisfy the equations:
\begin{equation}\begin{split}
a)\ & \frac{dK}{dt} = KH \\
b)\ & \frac{dH}{dt} = H^2-2K-4 \\
c)\ & \frac{dg}{dt} = -2gH.
\end{split}\end{equation}
\end{prop}

\noindent\underline{Proof of Proposition 4.1}: (a) and (b) follow immediately from (4.1) and (3.9). To prove (c) we use the standard result on ordinary differential equations: 

\noindent
\noindent\underline{Lemma} If A is an $n \times n$ matrix function of t and $a=\det\ A,$ then
$$\frac{da}{dt} = a\tr\left[A^{-1}\frac{dA}{dt}\right],$$
see \cite{C-L}.

Therefore:
\begin{eqnarray*}
\frac{dg}{dt} &=& g\tr\left[G^{-1}\frac{dG}{dt}\right] \\
&=& -2g\tr G^{-1}Q\\
&=& -2 g H.
\end{eqnarray*}
The last equality follows as $G^{-1}Q$ is the second fundamental form and thus:
$$\tr G^{-1}Q =  (k_1 +k_2)\ .$$
Using Cramer's rule to express $G^{-1}$ it is evident that $gG^{-1}$ is well
defined and differentiable even if $G$ is not invertible. \eproof

An important corollary of Proposition 4.1 is:

\begin{coro}  As long as $K(t)$ is finite
\begin{equation}
\frac{dK^2g}{dt}  =0.
\end{equation}
\end{coro}

\noindent\underline{Proof of Corollary 4.2}: We differentiate and substitute from (4.2) to obtain:
\begin{eqnarray*}
\frac{dK^2g}{dt} &=&  2K\frac{dK}{dt}g + K^2\frac{dg}{dt}  \\
&=& 2K^2Hg - 2K^2gH\\
&=& 0.
\end{eqnarray*}
\eproof

In fact $K^2(t)g(t)$ is constant whether or not $K$ becomes infinite. As we saw in $\S$3, $K(t)$ becomes infinite if and only if $g(t)$ goes to zero. An easy asymptotic analysis using formul{\ae} for $g_{ij}(t)$ and $K(t)$ following
from (3.12) and (3.10) respectively shows that the constant value of $K^2(t)g(t)$ does not change across a singularity of $K$.

Using formula (3.12), we get an even clearer picture of the behavior of the surface at a focal point. If a single principal curvature becomes infinite, then the corresponding principal vector $X_i(t)$ goes to zero. The orientation of the surface is therefore reversed as it passes through the focal point. If both curvatures become infinite, then the orientation is preserved. As these sign changes are mirrored by sign changes in $K$.

As $dA(t)=\sqrt{g(t)}dx\wedge dy,$ Corollary 4.2 can be restated:

\noindent\underline{Corollary 4.2'}: {\em $KdA$ is invariant under the parallel flow.}

This is important for it allows one to compute the Gauss curvature of $\Sigma(\rho)$ in terms of the first derivatives of $R_\rho$.
From the formul{\ae} (4.1) and (3.10) we derive expressions for $K(t)$ and $H(t)$. Let $K_0$ and $H_0$ denote the respective initial values.
\begin{equation}\begin{split}
K(t) & = \frac{K_0}{K_0\sh^2 t-H_0\sh t \ch t + (\ch^2t+\sh^2t)}\\
H(t) & = \frac{H_0(\sh^2 t +\ch^2t) -4 \sh t \ch t - 2K_0\sh t \ch t}{K_0\sh^2 t-H_0\sh t \ch t + (\ch^2t+\sh^2t)}.
\end{split}\end{equation}
From these formul{\ae} the following analogue of a classical theorem of Bonnet is immediate:

\begin{prop} If $\Sigma$ is a surface in $\Hsp$ with constant mean curvature $H$ such that
$$|H| > 2\ , $$
then some parallel surface $\Sigma_t$ has constant Gauss curvature $K> 0$ and vice versa.
\end{prop}
\noindent\underline{Proof of Proposition 4.3}
  The case that $|H_0|>2$ is immediate from the formula for $K(t),$ as there is always a $t$ such that
  \begin{equation}
    H_0=\frac{\ch t}{\sh t}+\frac{\sh t}{\ch t}.
  \end{equation}
  To prove the converse we differentiate the formula for $H(t)$ w.r.t. $H_0$ to find
  \begin{equation}
    \frac{\pa H(t)}{\pa H_0}=
      \frac{1-K_0\sh^2t}{[K_0\sh^2 t-H_0\sh t \ch t + (\ch^2t+\sh^2t)]^2}.
  \end{equation}
  If $K_0\sh^2t=1,$ then $H(t)$ does not depend on $H_0,$ which completes the proof.
\eproof

\noindent\underline{Remark}: The extreme values $|H| = 2$ and $K = 0$ correspond to surfaces with very simple defining equations in the representation described in $\S$3. See $\S$5.
\section{Explicit Formul{\ae} for Envelopes} In this section we return to the study of surfaces represented as
envelopes of horospheres. We will derive formul{\ae} for the first and
second fundamental forms of $\Sigma(\rho)$ and use these to study the behavior of $\Sigma_t(\rho)$ as $t \rightarrow \infty$.

Calculations will be done at a point $p$ on $\St$ about which geodesic normal coordinates have been introduced. At such a point:
\begin{eqnarray*}
1)\ g_{ij}^{\St}(p) &=& \delta^i_j\\
2)\ \Gamma^i_{jk}(p) &=& 0\ - \mbox{the Christoffel symbols of $g_{ij}^{\St}(p)$ }.
\end{eqnarray*}
In fact we will use coordinates provided by stereographic projection from the antipodal point to $p$, $p^*$. Such coordinates satisfy 1 and 2 at $p$. They also define local conformal parameters for $\St$. $p$ will always go to the point $z= 0$. The round metric on $\St$ is:
$$\frac{16|dz|^2}{(4+|z|^2)^2}\ .$$
By a rotation in $\R^3,$ $p$ can be normalized to be the point $(0,0,-1)$. The formula for $R_\rho(x,y) :U \rightarrow \Hsp$, $U$ an open set in $\C$, is:

\begin{equation}
\begin{split}
AR_\rho(x,y) &= \left[\frac{(4+r^2)^2}{16}(\rho_x^2+\rho_y^2)+e^{2\rho}-1\right]\left(\frac{4x}{4+r^2},\frac{4y}{4+r^2},\frac{r^2-4}{4+r^2}\right)\\
 & + 2\rho_x\left(1+\frac{y^2-x^2}{4},-\frac{xy}{2},x\right)+2\rho_y\left(-\frac{xy}{2}, 1+\frac{x^2-y^2}{4}, y\right) 
\end{split}
\end{equation}
where
\begin{eqnarray*}
r^2 &=& |z|^2\qquad\mbox{and}\\
A &=& \frac{(4+r^2)^2}{16}(\rho_x^2+\rho_y^2)+(e^{\rho}+1)^2\ .
\end{eqnarray*}
From (5.1) we can calculate the coordinate tangent vectors $\partial_xR_\rho$ and $\partial_yR_\rho$. Call them $X$ and $Y$ respectively. A very tedious but elementary calculation shows that:

\begin{equation}
\begin{split}
\left(
\begin{matrix} 
X\cdot X& X\cdot Y\\
X\cdot Y& Y\cdot Y
\end{matrix}\right)\Bigg|_{z=0} &=\\
\frac{1}{A^{2}}&\left(
\begin{matrix} 
2\rho_{xx}+\rho_y^2-\rho_x^2-1+e^{2\rho} & 2(\rho_{xy}-\rho_x\rho_y)\\
2(\rho_{xy}-\rho_x\rho_y)& 2\rho_{yy}+\rho_x^2-\rho_y^2-1+e^{2\rho}
\end{matrix}\right)^2
\end{split}
\end{equation}
at $z=0$ we have $A= (\rho_x^2+\rho_y^2)+(e^{\rho}+1)^2$.

This is the induced Euclidean metric for the immersion. The hyperbolic metric is easily obtained as
$$\langle S,T\rangle = 4S\cdot T(1- |R_\rho|^2)^{-2}\ .$$
From (5.1) we compute that:
\begin{equation}
\left.4(1- |R_\rho|^2)^{-2}\right|_{z=0} = \frac{A^2}{4e^{2\rho}}.
\end{equation}

Putting (5.2) and (5.3) together, we obtain:

\begin{prop}At the center, $p$ of a geodesic normal coordinate system on $\St$ the metric on $\Sigma\rightarrow \Hsp$ is given by:
\end{prop}
\begin{equation}
g_{ij}(p) =\left(\begin{matrix}
\frac{e^\rho}{2}+ (\rho_{xx} +\frac{\rho_y^2-\rho_x^2-1}{2})e^{-\rho}& (\rho_{xy}-\rho_x\rho_y)e^{-\rho}\\
(\rho_{xy}-\rho_x\rho_y)e^{-\rho}&\frac{e^\rho}{2}+ (\rho_{yy} +\frac{\rho_x^2-\rho_y^2-1}{2})e^{-\rho}
\end{matrix}\right)^2.
\end{equation}

To calculate the corresponding quantities for $\Sigma_t(\rho)$ it follows from Theorem 2.1 that one merely replaces $\rho$ by $\rho+t$ in (5.4). We write
\begin{equation}
g_{ij}(t,p) =\left(\begin{matrix}
\frac{e^{\rho+t}}{2}+Ee^{-(\rho+t)}& Fe^{-(\rho+t)}\\
Fe^{-(\rho+t)}&\frac{e^{\rho+t}}{2}+Ge^{-(\rho+t)}
\end{matrix}\right)^2
\end{equation}
where $E,F$ and $G$ are given by
\begin{eqnarray*}
E&=&\rho_{xx} +\frac{1}{2}(\rho_y^2-\rho_x^2-1)\\
F&=&\rho_{xy}-\rho_x\rho_y\\
G&=&\rho_{yy} +\frac{1}{2}(\rho_x^2-\rho_y^2-1).
\end{eqnarray*}
According to Theorem 3.1
$$\frac{dg_{ij}}{dt} = -2\Pi_{ij}\ .$$
Differentiating (5.5) and setting t= 0 we obtain:
\begin{equation}
\Pi_{ij} = \left(\begin{matrix}
\frac{e^{\rho}}{2}+Ee^{-\rho}& Fe^{-\rho}\\
Fe^{-\rho}&\frac{e^{\rho}}{2}+Ge^{-\rho}
\end{matrix}\right)\left(\begin{matrix}
Ee^{-\rho}-\frac{e^{\rho}}{2}& Fe^{-\rho}\\
Fe^{-\rho}&Ge^{-\rho}-\frac{e^{\rho}}{2}
\end{matrix}\right).
\end{equation}
The second fundamental form equals
$$\Pi^i_j = g^{ik}\Pi_{kj}.$$
Using (5.5) to compute $g^{ij}$ we obtain:

\begin{prop} At the center, $p$ of a geodesic normal coordinate
system on $\St$ the second fundamental form of $\Sigma(\rho)\rightarrow \Hsp$ is:
\begin{equation}
\Pi^i_j = D^{-1}\left[ \left(\begin{matrix}
\frac{E-G}{2}& F\\
F&\frac{G-E}{2}
\end{matrix}\right) +\left[(EG-F^2)e^{-2\rho}-\frac{e^{2\rho}}{2}\right] \Id
\right],
\end{equation}
here $D = (EG-F^2)e^{-2\rho}+\frac{e^{2\rho}}{4}+\frac{E+G}{2}$.
\end{prop}

\noindent\underline{Remarks}: \begin{enumerate} \item From formula (5.7) it is apparent that the principal directions are preserved under the parallel flow as the eigenvectors of $\Pi^i_j$ are clearly independent of $t$.
\item  From (5.5) it is clear that for a $C^2$-function $\rho$ the
$\det g_{ij}(t,p)$   is usually non-zero, In fact it follows
from our analysis of focal sets in $\S$3 that this determinant
vanishes at $0,1$ or $2$ values of $t.$
\item Therefore if $\rho \in C^2$, then $\Sigma(\rho)$ is either smooth at a
  given point $p$, or $\Sigma(\rho\pm t)$, for small $t$, is smooth at a point
  $q$ that projects to $p$ under the parallel flow. In the latter case $p$ is in
  the focal set of $\Sigma(\rho-t)$. From this we conclude that the only
  singularities that can arise in an envelope $\Sigma(\rho)$ generated by a
  $C^2$-function $\rho$, are those which arise as focal singularities of an
  immersed $C^1$-surface.
\end{enumerate}

Using (4.4) and (5.5), we can study the asymptotic behavior of $g_{ij}(t)$ and $K(t)$ as $t\rightarrow \infty$, A quick glance at (5.5) shows that $g_{ij}(t)$  tends uniformly to infinity  as $t\rightarrow \infty$. Thus we must rescale
to obtain a finite limit; let
\begin{equation}
\widetilde{g}_{ij}(t) = 4e^{-2t} g_{ij}(t).
\end{equation}
An elementary fact about Gauss curvature is that:
$$K(cg_{ij}) = c^{-1}K(g_{ij})$$
for c a positive constant.

Hence:
\begin{equation}
K(\widetilde{g}_{ij})= \frac{e^{2t}}{4} K(g_{ij}(t)).
\end{equation}
We will regard $\widetilde{g}_{ij}(t)$ as a family of metrics on a domain $\Omega \subset \St$. The Gauss curvatures of these metrics are calculated using (5.9), The asymptotic behavior is as follows:

\begin{prop}  Let $\rho$  be a  $C^4$-function on a domain $\Omega \subset \St$. Then
\begin{equation}\begin{split}
a)\qquad &  \lim_{t\rightarrow \infty} R_{\rho+t} = \Id,\\
b) \qquad & \lim_{t\rightarrow \infty} \widetilde{g}_{ij}(t) = e^{2\rho}d\sigma^2,\\
c) \qquad &\lim_{t\rightarrow \infty} K(\widetilde{g}_{ij}(t)) = \frac{K_0}{(K_0+2-H_0)}\\
  \qquad &= \frac{(k_1k_2-1)}{(1-k_1)(1-k_2)}.
\end{split}
\end{equation}
\end{prop}

\noindent\underline{Proof of Proposition 5.3}: The proof of part (a) follows from (2.4):

$$R_{\rho+t} = \frac{|D\rho|^2+(e^{2(\rho+t)}-1)}{|D\rho |^2+(e^{\rho+t}+1)^2}X(\theta) + \frac{2D\rho}{|D\rho|^2+(e^{\rho+t}+1)^2}.$$
As $\rho$ is differentiable
$$\lim_{t\rightarrow\infty} R_{\rho+t} = X$$
which is the assertion of part a.

At the center of a geodesic normal coordinate system on $\St$ ;  $\widetilde{g}_{ij}(t)$ is given by:

\begin{equation*}
\widetilde{g}_{ij}(t) = \left( \begin{matrix} 
e^\rho+2Ee^{-(\rho+t)} & 2Fe^{-(\rho+t)}\\
2Fe^{-(\rho+t)} & e^\rho+2Ge^{-(\rho+t)}\end{matrix} \right) ^2.
\end{equation*}

As $t\rightarrow \infty$ this tends to $e^{2\rho}\Id$. Since the metric is a tensor and the uniform limit of tensors is a tensor we have shown that the tensor
$$\lim_{t\rightarrow\infty} \widetilde{g}_{ij}(t)dx_i\otimes dx_j,$$
at the center of a geodesic normal coordinate system, is given by:
$$e^{2\rho} (dx_1\otimes dx_1 + dx_2\otimes dx_2)\  .$$
This tensor agrees at every point with the invariantly given tensor
$e^{2\rho}d\sigma^2$. This establishes (5.10)(b).  Using (5.9) and formula (4.4) we obtain:
$$K(\widetilde{g}_{ij}(t)) = \frac{K_0e^{2t}}{(K_0+2-H_0)e^{2t}}+O(e^{-2t}).$$

Letting $t \rightarrow\infty$ we obtain (5.10)(c). It is not immediate that the
limiting curvature function is the curvature of the limiting metric. We have
assumed that $\rho \in C^4(\Omega)$ and thus $\widetilde{g}_{ij}(t) \in
C^4(\Omega)$. From the form of $\widetilde{g}_{ij}(t)$ it is clear that the
second derivatives of the metric are equicontinuous in $t$. And therefore
$\widetilde{g}_{ij}(t)$ and its first two derivatives converge locally uniformly
to $e^{2\rho}d\sigma^2$ and its first two derivatives. Therefore the limiting
metric has a curvature which must coincide with the limiting value obtained
above. An elementary calculation shows that the two formul{\ae} given for this
limit agree.

\noindent\underline{Remarks:}
\begin{enumerate}
\item $\rho \in C^4$  is probably more restrictive than necessary as the formula for the curvature only involves the first two
derivatives of $\rho$. The functions we will be dealing with
are real analytic so we will not pursue the optimal smoothness hypothesis here.

\item Henceforth we will denote the limiting curvature by $K_\infty$ and the limiting area form by $dA_\infty$.

\item
 The virtue of formula (5.10)(c). is that the left hand side is a very simple expression in $\rho$ :
\begin{equation}
K_\infty = (1-\Delta_{\St})e^{-2\rho}.
\end{equation}
Whereas the right hand side is typically a fully non-linear
second order expression for most representations of immersed
surfaces. The other side of the coin is that the expressions
for \underline{both} $K$ and $H$ in our representation are fully non-linear second order quantities.

\item The three special classes of surfaces discussed in the introduction are given by the conditions:
\begin{eqnarray*}
a)\qquad K_\infty &=& 0,\\
b) \qquad K_\infty &=& 1,\\
c)\qquad  K_\infty &=& \alpha < 0.
\end{eqnarray*}

Putting these values into (5.10), we see that they correspond to an envelope that satisfies:
\begin{eqnarray*}
a')&  K &= 0,\\
b')&  H &= 2,\\
c')&  (1-\alpha)K &= \alpha(2-H) .
\end{eqnarray*}
Somewhat more general than (c') is the condition:
\begin{eqnarray*}
c'')\qquad (1-K_{\infty})K &=& K_{\infty}(2-H).
\end{eqnarray*}
for $K_\infty$, a negative function on $\Omega$. This gives rise to an
asymptotic Minkowski Problem for surfaces in $\Hsp$. In a
second paper we will study the regularity of the corresponding
surfaces. This study requires a rather detailed analysis
of the solution to (5.11) near $\partial\Omega$.
\end{enumerate}

We will study surfaces satisfying (c') in the last sections of this paper. Robert Bryant has written a beautiful paper on mean curvature 2 surfaces, [Br]. There he proves a result which we obtained independently.

\begin{prop} The Gauss map for a surface  is conformal if and only if $\Sigma$ has either mean curvature 2 or is umbilic.
\end{prop}

In our representation the proof of Proposition 5.4 is a tedious
calculation in local coordinates which we omit. We refer the interested reader to \cite{Br} where an elegant proof using moving frames is presented.

We will consider the following Dirichlet problem:
Given a collection of curves $\Gamma \subset \St$ which are the oriented boundary of a domain $\Omega \subset \St$ , find an immersed surface
$$\Sigma \rightarrow \Hsp$$
such that:

a) for an $\alpha<0,$ the principal curvatures of $\Sigma$ satisfy:
$$(k_1k_2-1) = \alpha(1-k_1) (1-k_2).$$
Note if $\alpha= -1$ this equation reduces to
$$k_1^{-1}+k_2^{-1}= 2\ ;$$
$\Sigma$ is a surface with mean radius of curvature $2$.


b ) $\overline\Sigma \cap\partial \Hsp = \Gamma$ ;  
$\overline\Sigma \cap\partial \Hsp$ is called the asymptotic boundary. It is denoted by $\partial_\infty \Sigma$.
$\Gamma$ must be oriented as the curvature equation involves the mean
curvature. We will think of $\Gamma$ as the boundary of a domain $\Omega$.  In
our representation this problem becomes a singular Dirichlet problem for a
complete conformal metric on $\Omega$ with curvature $\alpha$ :

Find a  $\rho \in C^2(\Omega)$ such that
\begin{equation}\begin{split}
a')& \qquad\ (1-\Delta_{\St}\rho)e^{-2\rho} = \alpha\\
b') & \qquad \lim_{p\rightarrow\partial\Omega} \rho = \infty.
\end{split}
\end{equation}
Observe that if we replace $\rho$ with $\rho_t=\rho+t,$ then
$$(1-\Delta_{\St}\rho_t)e^{-2\rho_t} = e^{-2t}\alpha.$$

If $\Sigma$ is smooth then its curvatures necessarily satisfy (a). From formula (2.4) it follows that (b') implies (b),

As a final result in this section we show that the invariance of the curvature form extends to $\infty$ :

\begin{prop} The curvature form $KdA$ equals the asymptotic curvature form $K_{\infty} dA_\infty$.
\end{prop}

\noindent\underline{Proof of Proposition 5.5}: The assertion follows immediately from Corollary 4.2' and the fact that
$$K(t)dA(t) = \widetilde{K}(t)d\widetilde{A}(t) \ .$$

The left hand side is computed with respect to the metric $g_{ij}(t)$ while the right hand side is computed with respect to $\widetilde{g}_{ij}(t)$. \eproof

The equation
\begin{equation}
KdA= K_{\infty} dA_\infty
\end{equation}
is an analogue of Gauss' formula for the curvature of a surface immersed in $\R^3$. Proposition 5.5 could be restated:

{\bf Proposition 5.5'}: {\em The pullback of the two form $K_\infty dA_\infty$, via the Gauss map of $\Sigma$ is the curvature form on $\Sigma$. }

\noindent\underline{Remarks}:
\begin{enumerate}
\item  The special case $H= 2 \Leftrightarrow K_\infty = 1$ is closest to the classical theorem, This case was also treated by Bryant.

\item (5.13) provides a relatively simple way to compute the Gauss
curvature of $\Sigma$ in terms of $\rho$. Using the formula for $K_\infty$, in terms of the principal curvatures, we see that:
\begin{equation}
dA = \frac{e^{2\rho}dA_{\St}}{(1-k_1)(1-k_2)}.
\end{equation}
\item  We close this section with the generating functions for several well known surfaces in $\Hsp$ :

a) A totally geodesic surface meeting $\St$ in an equator:
$\rho =  \log\frac{1}{\cos\theta},$
$\theta$ is the azimuthal angle relative to the normal direction
to the plane defining the equator.

b) A horosphere:
\begin{eqnarray*}
\rho &= &  \log\left(\frac{1}{\cos\theta}\right) +  \frac{1}{2}\log\left(\frac{1+\cos\theta}{1-\cos\theta}\right) +t\\
&=& \log\left(\frac{1}{1-\cos\theta}\right) + t,
\end{eqnarray*}

$\theta$ is the azimuthal angle measured from the point of tangency of the horosphere.
 
c) A geodesic connecting two antipodal points:

$$\rho =  \log\frac{1}{\sin\theta},$$

$\theta$ is the azimuthal angle measured from either endpoint
of the geodesic.
\end{enumerate}

\section{Conformal Parameters}

If we restrict ourselves to simply connected regions, $\Omega$ then the
Dirichlet problem in (5.12) is solved in terms of the conformal map from
$\Omega$ to the unit disk $D_1$ in $\C$. It is useful to express the fundamental
geometric quantities for $\Sigma(\rho)$ in terms of the complex derivatives of
$\rho$. As before, we use conformal parameters arising from stereographic
projection. When $z$ denotes such a parameter, the complex derivatives are
defined by:
\begin{equation}
\partial =\frac{1}{2}(\partial_x-i\partial_y)\qquad \overline\partial =\frac{1}{2}(\partial_x+i\partial_y).
\end{equation}
The round metric on the sphere is
\begin{equation}\begin{split}
d\sigma^2 &= \gamma^2|dz|^2\\
\gamma & = \frac{4}{4+|z|^2}.
\end{split}
\end{equation}

The expression for $g_{ij}$ in terms of complex derivatives is not
particularly illuminating. However, there is a Hermitian matrix with the same determinant and trace as $g_{ij}$ which is useful:

\begin{prop}
The matrix:
\begin{equation}
\left. h_{ij}\right|_{z=0} = e^{2\rho}\begin{bmatrix}
(2\dd \ddb \rho-1/2)e^{-2\rho}+1/2 &2(\dd^2 \rho-(\dd \rho)^2)e^{-2\rho}\\
2(\ddb^2 \rho-(\ddb \rho)^2)e^{-2\rho}&(2\dd \ddb \rho-1/2)e^{-2\rho}+1/2
\end{bmatrix}^2
\end{equation}
\end{prop}
has the same determinant and trace as $g_{ij}$ at $z=0$.

\noindent
The proof is an elementary calculation which we omit.

The formula for $\Pi^i_j$, (5.7) is considerably simpler in terms of the complex derivatives:

\begin{prop}
The second fundamental form is:

\begin{equation}
\left. \Pi^i_j\right|_{z=0} = D^{-1}\left[\left(\begin{matrix}
2\Re S & -2\Im S\\
-2\Im S & -2\Re S
\end{matrix}\right)+\lambda \Id\right]
\end{equation}
\begin{eqnarray*}
S &=& \dd^2\rho-(\dd\rho)^2\\
D &=& (EG-F^2)e^{-2\rho}+\frac{e^{2\rho}}{4}+\frac{E+G}{2}\\
\lambda &=& (EG-F^2)e^{-2\rho}-\frac{e^{2\rho}}{2}.
\end{eqnarray*}
\end{prop}

Again the proof is an elementary calculation which we omit.

\section{Weingarten Surfaces}

In this section we will solve the Dirichlet problem for $\rho$ and discuss the
regularity of the associated Weingarten surface, $\Sigma(\rho)$. We will obtain
formul{\ae} for the principal curvatures and identify the lines of curvature with
the trajectories of a holomorphic quadratic differential in $\Omega$. We will
suppose $\Omega$ is simply connected in this section,

Up to a scale factor, the Dirichlet problem, (5.12) is an equation for the logarithm of
the conformal factor for the complete hyperbolic metric on $\Omega \subset \St$.
Let $f_\Omega$ be a conformal map from $\Omega$ onto the unit disk, 
$D_1$ in $\C$. We obtain $\rho$ by pulling back the Poincar\'e metric on $D_1$ via $f_\Omega$. Let $z$ denote a conformal parameter on $\Omega$ and, differentiation with respect to $z$ then:
\begin{equation}
\rho = \log(2|f'_\Omega|\gamma^{-1}(1-|f_\Omega|^2)^{-1})\ .
\end{equation}
$\rho$  solves the equation:
$$(1-\Delta_{\St}\rho)e^{-2\rho} = -1\ .$$

Using Schwarz' Lemma and the Koebe $1/4$-theorem, we obtain the classical estimates, \cite{Ahl1}:
If $\Omega$ is a simply connected domain and $\partial\Omega$ has more than two points then
\begin{equation}
1/4\   \delta(z, \partial\Omega)^{-1} \leq e^\rho \leq  \delta(z, \partial\Omega)^{-1},
\end{equation}
$\delta(\cdot,\cdot)$ is Euclidean distance measured in the conformal parameter $z$. From
(7.2) it follows easily that there is an exhaustion of $\Omega$ by compact subregions $\Omega_n \subset \subset \Omega$ such that

\begin{equation}
\rho \geq  n  \mbox{ in } \Omega\backslash \Omega_n
\end{equation}
which implies that $\rho$ satisfies the boundary condition:
$$ \lim_{p\to \partial\Omega} \rho(p) = \infty\ .$$
This boundary condition is equivalent to
$$\overline{\Sigma(\rho)}\cap\St =\partial\Omega\ .$$

\begin{prop}
If $\rho \in C^2(\Omega)$ tends to $\infty$  on $\partial\Omega$ as in (7.3), then the surface $\Sigma(\rho)$ tends to $\St$ precisely along $\partial\Omega$.
\end{prop}

\noindent\underline{Remark}: For the proposition to be correct $\Sigma(\rho)$ does not need to be smooth.

\noindent\underline{Proof of Proposition 7.1}: Using (2.4), we calculate the Euclidean norms of $R_\rho$ and $R_\rho-X$ :
\begin{eqnarray*}
|R_\rho| &=&1 -\frac{4e^\rho}{|D\rho|^2+(e^\rho+1)^2}\\
|R_\rho-X|^2 &=& \frac{4(e^\rho-1)^2+4|D\rho|^2}{(|D\rho|^2+(e^\rho+1)^2)^2}.
\end{eqnarray*}

From the first formula it is clear that $R_\rho(z)$ lies in the interior of $\Hsp$ for $z$ in the interior of $\Omega$. From the second formula we easily obtain that
\begin{eqnarray*}
|R_\rho-X| &\leq&  2 \min (|D\rho|^{-1}, |D\rho|(e^\rho+1)^{-2})\\
&\leq& 2(e^\rho+1)^{-1}.
\end{eqnarray*}
Hence (7.3) implies:
\begin{equation}
|R_\rho-X| \leq 2e^{-n} \mbox{ for } n\in  \Omega\backslash \Omega_n.
\end{equation}

This completes the proof of the proposition. \eproof

To study the regularity and lines of curvature of $\Sigma(\rho)$ , we rewrite
the formul{\ae} for $\Pi^i_j$ and $h_{ij}$ in terms of $f_\Omega$.

\begin{prop}
 At the center of a stereographic coordinate system on $\St$, $(z=0)$ $h_{ij}(t)$  and $\Pi^i_j(t)$ are given by:
\begin{equation}\begin{split}
h_{ij}(t) &= e^{2\rho+t}\left[\begin{matrix}
1/2(1+e^{-2t})& \Sh_{f_\Omega}(0)\mu_\Omega(0)^{-1}e^{-2t}\\
\overline{\Sh_{f_\Omega}(0)}\mu_\Omega(0)^{-1}e^{-2t}&1/2(1+e^{-2t})\end{matrix}\right]^2 \\
&\Pi^i_j(t) = D^{-1}\left[\left(\begin{matrix}
\Re \Sh_{f_\Omega}(0)& -\Im \Sh_{f_\Omega}(0)\\
-\Im \Sh_{f_\Omega}(0)&-\Re \Sh_{f_\Omega}(0)\end{matrix}\right)+\lambda \Id\right].
\end{split}
\end{equation}
\end{prop}

Here  $\Sh_f(z) =(\frac{f''}{f'})' - \frac{1}{2}(\frac{f''}{f'})^2$ is the Schwarzian derivative of $f$ with respect to the conformal parameter $z$; the conformal factor for the hyperbolic metric is:
$$\mu_\Omega = 4|f'_\Omega|^2(1-|f_\Omega|^2)^{-2}.$$

\noindent\underline{Remark}: The function of $\Sh_f(z)$ is not well defined in $\Omega$; however, the quadratic differential:
$$q_\Omega = \Sh_{f_\Omega}(z)dz^2$$
is well defined. For the case of Weingarten surfaces, the Gauss map is not conformal and therefore $q_\Omega$ does not define a holomorphic differential on $\Sigma(\rho)$.

\noindent\underline{Proof of Proposition 7.2}: To derive (7.5) from (6.3) and (6.4), we need only calculate:

$$\left.(2\dd \ddb\rho - 1/2)e^{-2\rho}\right|_{z=0}$$
and
$$\left.(\dd^2\rho -(\dd \rho)^2)e^{-2\rho}\right|_{z=0}$$
in terms of $f_\Omega$. Using the facts that
$$\left.\gamma\right|_{z=0} = 1$$
and
$$\left.\dd\gamma\right|_{z=0} = \left.\dd^2\gamma\right|_{z=0} =0$$
we easily obtain:
$$(2\dd\rho \ddb\rho - 1/2)e^{-2\rho} = 1/2$$
and
$$\left.e^{2\rho}\right|_{z=0} = \left. 4|f'_\Omega|^2(1-|f_\Omega|^2)^{-2}\right|_{z=0}\ .$$
A straightforward calculation shows:
$$\left.(\dd^2 \rho-(\dd \rho)^2)\right|_{z=0} = \Sh_{f_\Omega}(0)\ .$$
\eproof

\noindent\underline{Remark:} The last formula appears in \cite{Br}.

As $h_{ij}(t)$ is a matrix of the form:
$$e^{2\rho_t}[Ae^{-2t}+1/2],$$
a bound on the eigenvalues of $A,$ uniform in $\Omega,$ would imply that 
$$\det h_{ij}(t) > 0$$
for t sufficiently large. As $h_{ij}$ and $g_{ij}$ have the same determinant, it would then follow that $R_{\rho+t}$ is an immersion for sufficiently large $t$.

The estimate we need follows from a theorem of Kraus; \cite{Kr}:

\noindent\underline{Theorem A}: {\em If $f$ is a univalent function in the disk, then
\begin{equation}
|\Sh_{f}(z)| \leq 6(1-|z|^2)^{-2}.
\end{equation}
}

From this theorem one easily deduces:

\begin{prop} If $f_\Omega$ is a conformal map from a simply connected domain $\Omega$ onto the unit disk, then:
\begin{equation*}
|\Sh_{f_\Omega}(z)| \leq \frac{3}{2}\mu_{\Omega}(z).\qquad\qquad (7.6')
\end{equation*}
\end{prop}

\noindent\underline{Proof of Proposition 7.3}: By expressing $\Sh_f$ in terms of $\Sh_{f^{-1}}$ and recognizing the right hand side as $3/2\mu_{D_1}$,  (7.6') follows from (7.6). \eproof

To apply this estimate we need an expression for the principal curvatures in terms of $\Sh_{f_\Omega}$ and $\mu_\Omega$. For the remainder of this section we will use the notation $\Sh_\Omega$ for $\Sh_{f_\Omega}$.

\begin{prop}
 The principal curvatures of the $-1$-Weingarten surface, $\Sigma(\rho)$ are: 
 \begin{equation}
   k_\pm = \frac{1}{1\pm \left|\frac{\mu_\Omega}{\Sh_\Omega}\right|}=\frac{\frac{|\Sh_\Omega|}{\mu_\Omega}}
   {\frac{|\Sh_\Omega|}{\mu_\Omega}\pm 1}.
\end{equation}
\end{prop}


\noindent\underline{Proof of Proposition 7.4}: We begin with the formula derived in \S 5:
$$KdA = K_\infty dA_\infty\ .$$
Recall that $K_\infty=-1$ and $dA_\infty =\mu_\Omega dxdy$. On the other hand, from Proposition 7.2 it follows that:
\begin{equation}\begin{split}
\left.dA\right|_{z=0} &= \sqrt{\det h_{ij}} dxdy\\
&= \mu_\Omega\left(1-\left|\frac{\Sh_\Omega}{\mu_\Omega}\right|^2\right)dxdy.
\end{split}
\end{equation}
At $z=0$
$$K = \left(\left|\frac{\Sh_\Omega}{\mu_\Omega}\right|^2 - 1\right)^{-1}\ .$$
K is given terms of the principal curvatures by
\begin{equation}
K = k_1k_2-1\ .
\end{equation}
As $\Sigma(\rho)$ is $-1$-Weingarten surface its principal curvatures satisfy the relation:
\begin{equation*}
\frac{1-k_1k_2}{(1-k_1)(1-k_2)} = 1
\end{equation*}
or
\begin{equation}
k_2= \frac{k_1}{2k_1-1} \ .
\end{equation}
Substituting into (7.9) from (7.10), we easily obtain:
\begin{equation}
k_\pm = K+1\pm\sqrt{K(K+1)}\  .
\end{equation}
(7.7) follows from (7.11) after substitution from (7.8). 
The formula was derived for $z=0$. However, $\left|\frac{\Sh_\Omega}{\mu_\Omega}\right|$ is invariant under conformal changes of parameter and so the formula for $K$ is valid throughout the domain of the conformal parameter.

The regularity theorem for Weingarten surfaces over simply connected regions is:

\begin{thm}
 Suppose $\Omega \subset \St$ is a simply connected region,
\begin{enumerate}[label=\alph*)] 

\item  If  $|\Sh_\Omega| < \frac 12\mu_\Omega$ everywhere in $\Omega$ then the $\alpha$-Weingarten
surfaces with boundary equal to $\partial\Omega$ are smoothly embedded for all $\alpha < 0$.
\item If  $|\Sh_\Omega| < \mu_\Omega$ everywhere in $\Omega$  then the $\alpha$-Weingarten
surfaces with boundary equal to $\partial\Omega$ are immersed for $-1<\alpha <0$.

\item If $\Omega$ satisfies neither (a) nor (b), then the $\alpha$-Weingarten surfaces with boundary equal to $\partial\Omega$ are immersed for:
$$-\frac{1}{2}<\alpha<0.$$
\end{enumerate}
\end{thm}

\noindent\underline{Proof of Theorem 7.5}: We use the fact that if
$\Sigma(\rho)$ is a $-1$-Weingarten surface, then $\Sigma(\rho_t)$ is a
$-e^{-2t}$-Weingarten surface.  The proof follows from the formula for the
hyperbolic Jacobian determinant of $R_{\rho+t},$
$$dA(t) = -dA_\infty K(t)^{-1},$$
and the formula for the curvature:
$$K(t) = \left(\left|\frac{\Sh_\Omega}{\mu_\Omega}\right|^2 e^{-2t} - \ch^2 t\right)^{-1}\ .$$
Since $dA_\infty$, is always positive and $K(t)$ is negative and finite whenever
\begin{equation}
|\Sh_\Omega |^2  < \mu^2_\Omega e^{2t}\ch^2 t\ ,
\end{equation}
assertion (b) follows immediately from (7.12). Assertion (c) follows from (7.12) and Proposition (7.3).

To prove assertion (a) we use the formula for the principal curvatures of a $-1$-Weingarten surface, (7.7):
$$k_\pm =
\frac{\left|\frac{\Sh_\Omega}{\mu_\Omega}\right|}{\left|\frac{\Sh_\Omega}{\mu_\Omega}\right|\pm
  1}.$$ If $|\Sh_\Omega | < 1/2 |\mu_\Omega |$ then $|k_\pm | < 1$, hence
$\Sigma(\rho)$ is a smooth immersion of $\Omega$ without boundary points in
$\Hsp$ and thus complete. Theorem 3.4 applies to show that the $-1$-Weingarten
surface is embedded as well as the family of surfaces parallel to it.  \eproof

\noindent\underline{Remark}: If $\Omega$ satisfies neither (a) nor (b), then the $-1$-Weingarten surface with boundary equal to $\partial\Omega$ has singularities. This could only fail to occur if $|\Sh_\Omega | \mu_\Omega ^{-1} > 1$ everywhere in $\Omega$. If we choose our conformal parameter by stereographically projecting from a point in the interior of $\Omega$, then $\partial\Omega$ will be a compact set and $\mu_\Omega$ will be bounded from below near $\partial\Omega$. Thus $\Sh_\Omega(z)^{-1}$ will be holomorphic near $\partial\Omega$ and satisfy the estimate:
$$|\Sh_\Omega(z)^{-1}| \leq \delta(z,\partial\Omega)\ ,$$
clearly an absurdity for a holomorphic function.

In case (a) of the theorem the parallel surfaces $\Sigma_t$ define a foliation of a part of $\Hsp$. To prove this, we need a comparison theorem for surfaces:

\begin{prop}
 Let $\Omega_1$ and $\Omega_2$ be two simply connected domains on $\St$ and suppose
$$\Omega_1 \subset\subset\Omega_2\ .$$
Let $\rho_1$ and $\rho_2$ define complete metrics on $\Omega_1$ and $\Omega_2$, respectively such that:
$$\rho_2 < \rho_1 \mbox{ in } \overline{\Omega}_1\ .$$
Finally, suppose that the surfaces $\Sigma(\rho_2)$ and $\{\Sigma_t(\rho_1)\  :\ t > 0\}$ are properly embedded, then
$$\Sigma_t(\rho_1)\cap \Sigma(\rho_2)=\emptyset\  .$$
\end{prop}

\noindent\underline{Proof of Proposition 7.6}: First we see that there exists a $t_0 > 0$  such that
$$\Sigma_t(\rho_1)\cap \Sigma(\rho_2)=\emptyset\mbox{ if } t\geq t_0\ .$$
If this were not the case, we could find a sequence of times $t_n \rightarrow \infty$ and a sequence of points:
$$p_n\in \Sigma_{t_n}(\rho_1)\cap \Sigma(\rho_2)\ .$$
As $\B^3$ is compact, there is a subsequence of the points $p_n$ converging to $p^*$; call this sequence $p_n$ well. From formula (2.4), it is evident that $p^*$ lies in $\overline{\Omega}_1$. On the other hand, $\Sigma(\rho_2)\cap\partial \B^3 = \partial\Omega_2$, thus we have derived a contradiction as
$$\Omega_1 \subset\subset\Omega_2\ $$
and $p^* \in \Sigma(\rho_2)\cap\partial \B^3$.

As $\Sigma(\rho_2)$ is properly embedded it divides $\B^3$ into two connected components, $D_1$ and $D_2$. We label them so that $\Omega_1$ lies in $\overline{D}_1$ . From the above argument, it follows that
\begin{equation}
\Sigma_t(\rho_1)\subset\subset D_1
\end{equation}
for $t$ large enough. As $\partial\Omega_1$ is disjoint from $\partial\Omega_2$ , there is a first time $t_1$ such that
$$\Sigma_{t_1}(\rho_1)\cap \Sigma(\rho_2)\neq\emptyset\ .$$
These two surfaces are tangent at some finite point, $q$. Since both $\Sigma(\rho_2)$ is embedded and (7.13) holds, it follows that the inward pointing
normals of $\Sigma_{t_1}(\rho_1)$ and $\Sigma(\rho_2)$ agree at $q$. Thus the Gauss maps of the two surfaces agree at $q$ ; call the common value $\theta$. From the definitions of these surfaces as envelopes of horospheres, it follows that
$$H(\theta,\rho_1(\theta)+t_1) = H(\theta,\rho_2(\theta))\ .$$
From this we conclude that
\begin{eqnarray*}
t_1 &=& \rho_2(\theta)-\rho_1(\theta)\\
&< & 0\ .
\end{eqnarray*}
Therefore, $\Sigma(\rho_1)\cap \Sigma(\rho_2)= \emptyset$. \eproof

\begin{coro} The $\alpha$-Weingarten surfaces for $\alpha < 0$ and $\Omega$ satisfying assertion (a) of Theorem 7.5 foliate  $\bigcup_{t\in\R} \Sigma_t(\rho)$.
\end{coro}

\noindent\underline{Proof of Corollary 7.7}: Proposition 7.6 does not apply immediately as all the surfaces $\Sigma_t(\rho)$ share a common boundary, To remedy this, we consider the subregion $\Omega_r$ of $\Omega$ given by:
$$\Omega_r  = f_\Omega^{-1}\{ z\ : |z| < r \}.$$
The conformal map for $\Omega_r$ onto the unit disk is
$$f_r = \left. r^{-1}f_{\Omega}\right|_{\Omega_r}$$
An easy calculation shows:
$$\Sh_r = \Sh_{f_r}= \Sh_\Omega$$
and the comparison principle for hyperbolic metrics [Ahl2] states:
$$\mu_\Omega < \mu_{\Omega_r} \mbox{ on } \Omega_r\ $$
Thus:
\begin{eqnarray*}
\left|\frac{\Sh_r}{\mu_{\Omega_r}}\right| &=&\left|\frac{\Sh_\Omega}{\mu_{\Omega_r}}\right|\\
&<&\left|\frac{\Sh_\Omega}{\mu_{\Omega}}\right|\\
&\leq&\alpha,
\end{eqnarray*}
for $\alpha = \sup_{\Omega}\left|\frac{\Sh_\Omega}{\mu_{\Omega}}\right| < 1/2$. Thus, for all $r\leq 1$:
$$\left|\frac{\Sh_r}{\mu_{\Omega_r}}\right| < 1/2 \mbox{ in } \Omega_r.$$

Letting $2\rho_r = \log \mu_{\Omega_rr}\gamma^{-2}$ it follows from Theorem 7.5 (a) that $\Sigma_t(\rho_r)$  is embedded for every $t$. Proposition 7.6 applies and we conclude that $\Sigma_t(\rho_r)$  and $\Sigma_s(\rho)$ are disjoint for $t> s$ and $r$ less than $1$.
Letting $r\rightarrow 1$, it follows that $\Sigma_t(\rho_r)\cap \Sigma_s(\rho)$ consists of points of
tangency; at the common image of such points under the respective Gauss maps, $\theta$:
$$\rho(\theta) +t = \rho(\theta) + s,$$
an obvious contradiction to the choice of $t$ and $s$.
\eproof

\noindent\underline{Remarks}: 
\begin{enumerate}
\item For general regions $\Omega$ with $|\Sh_\Omega|\mu_\Omega^{-1} < 1/2$ we do not yet know if $\bigcup_t \Sigma_t$  is all of $\Hsp$. This will be the case if $|\nabla\rho|\rightarrow \infty$ on $\partial\Omega$. The map $F(\theta)= \lim_{t\rightarrow -\infty}  R_{\rho+t}(\theta)$ would then define a homeomorphism of $\Omega$ onto $\Omega^c$ which fixes $\partial\Omega$ pointwise. A simple degree argument then shows that $\bigcup_t \Sigma_t= \Hsp$. It seems likely that this is true in full generality, but we have not yet obtained uniform lower estimates for $|\nabla\rho|$ near $\partial\Omega$.

\item The estimate $|\Sh_\Omega|\mu_\Omega^{-1} < 1/2$ holds whenever $\Omega$
  is the stereographic image of a convex region in $\C$ , [Ne]. Estimates of the
  form $|\Sh_\Omega|\mu_\Omega^{-1} < m$ for an $m< 3/2$ follow if $f$ has a
  quasi-conformal extension to all of $\St$ ; \cite{AhlWe} and \cite{Pom}.
\end{enumerate}
The final result in this section describes the lines of curvature on $\alpha$-Weingarten surfaces.

\begin{thm} The lines of curvature on an $\alpha$-Weingarten surface project under the Gauss map to the positive and negative trajectories of the holomorphic quadratic differential; $q_\Omega = \Sh_\Omega(z) dz^2$.
\end{thm}

\noindent\underline{Remark}: This theorem holds for all values of $\alpha$. For $\alpha > 0$ a function $\rho$ defining the surface is locally expressible in terms of a holomorphic function $f$ and  $\Sh_\Omega = \Sh_f$ locally. We will only treat the case $\alpha < 0$. 

Note $q_\Omega$ is holomorphic on $\Sigma(\rho)$ if and only if $\alpha=1$.

\noindent
\noindent\underline{Proof of Theorem 7.8}: The positive and negative trajectories are the integral curves of the line fields defined by:
$$\Im q_\Omega= 0\ .$$
If $\Sh_\Omega(z) = a+ib$ then line fields at $z$ are given by:
$(a - \sqrt{a^2 +b^2} , -b)$ and $(a +\sqrt{a^2 +b^2} , -b)$. The principal directions are determined at $z=0$ by the matrix:\\
$$
\left(\begin{matrix}
\Re\Sh_\Omega& -\Im\Sh_\Omega\\
-\Im\Sh_\Omega&-\Re\Sh_\Omega
\end{matrix}\right) = 
\left(\begin{matrix}
a& -b\\
-b&-a
\end{matrix}\right).
$$

The eigenvectors are easily seen to coincide with the vectors determining the trajectories of $q_\Omega$.
\eproof

\section{Higher Connectivity}
In the previous section, the solution of the Dirichlet problem for a Weingarten
surface over a simply connected region, $\Omega$ is reduced to the construction
of a complete hyperbolic metric on $\Omega$. The regularity depends upon an
upper bound for $|\Sh_\Omega|\mu_\Omega^{-1}$.

An essentially arbitrary domain, $\Omega$ on $\St$ has a complete hyperbolic metric, This is a consequence of the general uniformization theorem, If $\St\backslash \Omega$ consists of more than two points, then there is a conformal covering map $g$ from the unit disk to $\Omega$. $g$ is locally univalent so we can
define a local inverse $f:\Omega\rightarrow D_1$. If $f_1$ is a different local inverse, then there exists $\alpha$ and $\beta$ such that:
$$f_1 = (\alpha f-\beta)(\overline\beta f- \overline\alpha)^{-1}$$
$$|\alpha|^2 -|\beta|^2 =1.$$
Thus the hyperbolic metric for $\Omega$ can be expressed without ambiguity by:
$$e^{2\rho}d\sigma^2 = \frac{4|f'|^2}{(1-|f|^2)^2}|dz|^2.$$
The quadratic differential,
$$q_\Omega(z) = \Sh_f(z) dz^2,$$
is also well defined because the Schwarzian derivative is invariant under M\"obius transformations. Moreover, the formula for the curvature:
$$K = \left(\left| \frac{\Sh_f}{\mu_\Omega}\right|^2-1\right)^{-1}$$ is still
valid. Thus the smoothness of the $\alpha$-Weingarten surfaces over $\Omega$ is
still determined by the ratio $|\Sh_f|\mu_\Omega^{-1}$.

In [Ge] the following lemma is proved for arbitrary domains $\Omega$:

\begin{lemma} If $f$ is univalent in $\Omega$ then
\begin{equation}
|\Sh_f(z)| \leq 6\delta(z,\partial\Omega)^{-2}.
\end{equation}
\end{lemma}
Recall $\delta(\cdot,\cdot)$ is Euclidean distance in the conformal parameter
$z$. As the proof only requires an estimate of $\Sh_f$ in a disk contained in
$\Omega$, it is clear that $f$ need not be single valued. Gehring's argument
actually proves:

\noindent\underline{Lemma 8.1}: {\em If $f$ is holomorphic and locally univalent though not necessarily single valued in $\Omega$ then:
$$|\Sh_f(z)| \leq 6\delta(z,\partial\Omega)^{-2}\ . \qquad (8.1)$$
}

To prove regularity for the $\alpha$-Weingarten surfaces, we need a lower bound for $\mu_\Omega$ in terms of 
$\delta(z,\partial\Omega)$. For a general domain no such estimate is true, However, if $\Omega$ has the property that
$$\partial \Omega = \bigcup \gamma_i$$
where each $\gamma_i$ has more than two points and the component of $\St\backslash \gamma_i$ which
contains $\Omega$ is simply connected, then such an estimate holds. Let the component of $\St\backslash \gamma_i$, described above be denoted by $\Omega_i$. Estimate (7.2) applies to
the domains $\Omega_i$ to give:
\begin{equation}
1/4\ \delta(z,\partial\Omega_i)^{-2} \leq \mu_{\Omega_i}.
\end{equation}

We use the comparison principle for hyperbolic metrics to conclude:

\begin{prop} Let $\Omega$ be as described above, then
\begin{equation}
1/4\ \delta(z,\partial\Omega)^{-2} \leq \mu_{\Omega}.
\end{equation}
\end{prop}

\noindent\underline{Proof of Proposition 8.2}: The comparison principle for
hyperbolic metrics [Ahl2] states that if $\Omega_1\subset \Omega_2$ then
$$\mu_{\Omega_1} \geq \mu_{\Omega_2}\ .$$
As $\Omega \subset \Omega_i$ for every $i$ and therefore
$$\mu_{\Omega_i} \leq \mu_{\Omega} \mbox{ for every i } \ .$$
We apply (8.2) and take the supremum over $i$ to obtain:
$$\frac{1}{4}\sup_i \delta(z,\partial\Omega_i)^{-2} \leq \mu_\Omega(z).$$
The estimate, (8.3) follows as $\bigcup \gamma_i = \partial\Omega$. \eproof

\begin{coro}
For $\Omega$ as described above:
\begin{equation}
|\Sh_\Omega|\mu_\Omega^{-1} \leq 24.
\end{equation}
\end{coro}

As a corollary of the corollary, we have:

\begin{coro} If $\Omega$ is as described above, then the $\alpha$-Weingarten surfaces with boundary equal to $\partial\Omega$ are smoothly immersed if
$$-\frac{1}{47}<\alpha <0.\ $$
\end{coro}
The argument is identical to that used in the proof of  Theorem 7.5(c).

\end{document}